\documentclass{amsart}
\usepackage{amsfonts}
\usepackage{amssymb}
\usepackage{amsmath}
\usepackage[dvips]{epsfig,graphics}
\usepackage{latexsym, amsfonts, longtable, tabularx}
\usepackage{epsfig}
\usepackage{psfrag}
\usepackage{amscd}
\usepackage[all]{xy}

\newtheorem{thm}{Theorem}
\newtheorem{prop}{Proposition}
\newtheorem{lemma}{Lemma}

\newtheorem{defn}{Definition}
\newtheorem{rmk}{Remark}

\newcommand{\Z}{\mathbb{Z}}

\newcommand{\T}{\mathcal{T}}
\renewcommand{\c}{\cite}

\newcommand{\ZA}{{\mathbb{Z}}[A^{\pm1}] }
\newcommand{\pf}{{\em Proof: \quad }}
\newcommand{\done}{\hfill $\blacksquare$}
\newcommand{\kb}[1]{\ensuremath{\langle #1 \rangle}}

\newcommand{\del}{\partial}
\newcommand{\U}{{\widetilde U}}

\newcommand{\UC}{\mathcal{UC}}
\newcommand{\C}{\mathcal{C}}

\begin{document}

\title{Spanning trees and Khovanov homology}
\author{Abhijit Champanerkar}
\address{Department of Mathematics and Statistics, University of South Alabama}
\email{achampanerkar@jaguar1.usouthal.edu}
\thanks{The first author is supported by NSF grant DMS-0455978}

\author{Ilya Kofman}
\address{Department of Mathematics, College of Staten Island, City University of New York}
\email{ikofman@math.csi.cuny.edu}
\thanks{The second author is supported by grants NSF DMS-0456227 and PSC-CUNY 60046-3637}

\date{May 20, 2007}

\begin{abstract}
\noindent
The Jones polynomial can be expressed in terms of spanning trees of
the graph obtained by checkerboard coloring a knot diagram.  We show
there exists a complex generated by these spanning trees whose
homology is the reduced Khovanov homology.
The spanning trees provide a filtration on the reduced Khovanov
complex and a spectral sequence that converges to its homology.
For alternating links, all differentials on the spanning tree complex are zero and
the reduced Khovanov homology is determined by the Jones polynomial
and signature. We prove some analogous theorems for (unreduced) Khovanov homology.
\end{abstract}
\maketitle

\section{Introduction}

For any diagram of an oriented link $L$, Khovanov \c{Khovanov}
constructed bigraded abelian groups $H^{i,j}(L),$ whose
bigraded Euler characteristic gives the Jones polynomial $V_L(t)$:
$$\chi(H^{i,j})= \sum_{i,j}(-1)^i q^j {\rm rank}(H^{i,j})= (q+q^{-1})V_L(q^2) $$
Khovanov's homology groups turn out to be stronger invariants than the Jones polynomial.
For knots, Khovanov also defined a reduced homology
$\widetilde{H}^{i,j}(L)$ whose bigraded Euler characteristic is
$q^{-1}V_L(q^2)$ \c{KhPatterns}.

The Jones polynomial has an expansion in terms of spanning trees of
the Tait graph, obtained by checkerboard coloring a given link diagram
\c{thistlethwaite}.  Every spanning tree contributes a monomial to
the Jones polynomial, and for alternating knots, the number of
spanning trees is exactly the $L^1$-norm of the coefficients of the Jones
polynomial. These monomials are Kauffman brackets of certain twisted unknots (Theorem \ref{sum_unknots}).

We show the reduced Khovanov complex $\widetilde{C}(D)$ retracts to a {\em spanning
tree complex,} whose generators correspond to spanning trees of the Tait graph (Theorem \ref{mainthm}).
The main idea is that a spanning tree corresponds to a twisted unknot $U$,
and $\widetilde{C}(U)$ is contractible, providing a deformation retract of $\widetilde{C}(D)$.
This extends to (unreduced) Khovanov homology (Theorem \ref{unred}). The proof
does not provide an intrinsic description of the differential on spanning
trees. From a partial order on spanning trees, we get a filtration of $\widetilde{C}(D)$,
and a spectral sequence that converges to $\widetilde{H}(D)$
(Theorem \ref{SS}).

A knot $K$ is alternating if and only if all the spanning trees are in
one row of the spanning tree complex and hence all differentials on the spanning tree complex are zero.
We give a simple new proof that for alternating links $\widetilde{H}(K)$
is determined by its Jones polynomial
and signature (Theorem \ref{altredKH}).  We also give simple new
proofs for theorems in \c{Lee, AP, Manturov}
on the support of Khovanov homology of alternating and $k$-almost
alternating knots (Theorem \ref{thick}).

Wehrli independently gave a spanning tree model for Khovanov homology in \c{Wehrli}.

{\bf Acknowledgments}\quad
We are grateful to Oleg Viro for many contributions. We thank the
organizers of {\em Knots in Poland 2003}, and we also thank Alexander
Shumakovitch and Peter Ozsv\'ath for useful discussions.

\section{Spanning trees and twisted unknots}

There is a $1$-$1$ correspondence between connected link diagrams $D$
and connected planar graphs $G$ with signed edges.  $G$ is obtained by
checkerboard coloring complementary regions of $D$, assigning a vertex
to every shaded region, an edge to every crossing and a $\pm$ sign to
every edge such that for a positive edge, the $A$-smoothing joins the
shaded regions.  The signs are all equal if and only if $D$ is
alternating.  $G$ is called the Tait graph of $D$. Thistlethwaite
\cite{thistlethwaite} described the following expansion of the Jones
polynomial in terms of spanning trees of the Tait graph.

Fix an order on edges of $G$.  For every spanning tree $T$
of $G$, each edge $e$ of $G$ has an activity with respect to $T$,
as originally defined by Tutte.
If $e \in T$, $\mathit{cut(T,e)}$ is the set of edges that connect $T\setminus e$.
If $f \notin T$, $\mathit{cyc(T,f)}$ is the set of edges in the unique cycle of $T \cup f$.
Note $f \in cut(T,e)$ if and only if $e \in cyc(T,f)$.
An edge $e \in T$ is called internally active with respect to $T$ if
it is the lowest edge in its cut, and otherwise it is internally
inactive.  An edge $e \notin T$ is externally active with
respect to $T$ if it is the lowest edge in its cycle, and otherwise
it is externally inactive.  Each edge
has one of eight possible activities, depending on whether $(i)\; e
\in T$ or $e \notin T$, $(ii)\; e$ is active (live) or inactive (dead), $(iii)\; e$
has $\pm$ sign.
Let $L,\ D,\ \ell,\ d$ denote a positive edge that is internally
active, internally inactive, externally active,
externally inactive, respectively.  Let $\bar{L},\ \bar{D},\
\bar{\ell},\ \bar{d}$ denote activities for a
negative edge.
Each edge $e$ of $G$ is assigned a monomial $\mu_e\in\ZA$, as in Table \ref{Table2} (last row).
Let $\mu(T)=\prod_{e\in G}\mu_e$.

\begin{thm}[\cite{thistlethwaite}]\label{ThistThm}
Let $D$ be any connected link diagram.  Let $G$ be its Tait graph with
any order on its edges. Then the Kauffman bracket
$\kb{D}=\sum_{T \subset G} \mu(T)$.
\end{thm}

\begin{table}[h]
\caption{Activity word for a spanning tree determines a twisted unknot}
\label{Table2}
\begin{center}
\begin{tabular}{cc|cc|cc|cc}
$L$ & $D$ & $\ell$ & $d$ & $\bar{L}$ & $\bar{D}$ &
$\bar{\ell}$ & $\bar{d}$ \\
\hline
$-$ & $A$ & $+$ & $B$ & $+$ & $B$ & $-$ & $A$ \\
\includegraphics[height=0.5cm]{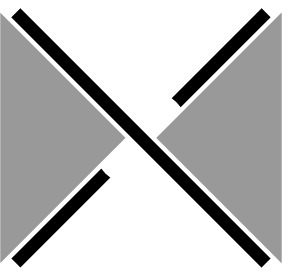} & \includegraphics[height=0.5cm]{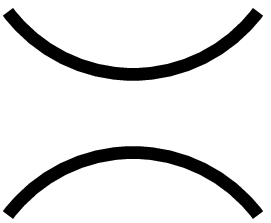} &
\includegraphics[height=0.5cm]{poscross.eps} & \includegraphics[height=0.5cm]{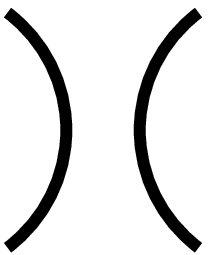} &
\includegraphics[height=0.5cm]{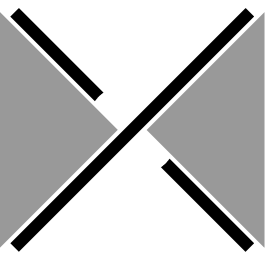} & \includegraphics[height=0.5cm]{posA.eps} &
\includegraphics[height=0.5cm]{negcross.eps} & \includegraphics[height=0.5cm]{posB.eps} \\
\hline
$-A^{-3}$ & $A$ & $-A^3$ & $A^{-1}$ & $ -A^{3}$ & $A^{-1}$ &
$-A^{-3}$ & $A$ \\
\end{tabular}
\end{center}
\end{table}

Using the edge order, we write an {\em activity word} for each
spanning tree $T$ using the eight letters for its edge activities.
$T$ is given by the capital letters of the word.
A \emph{twisted unknot} $U$ is obtained from the round unknot using only Reidemeister I moves.

\begin{lemma}\label{treetounknot}
Given an activity word for a spanning tree $T$, changing the crossings of
$D$ according to Table \ref{Table2} for dead edges and leaving the crossings
unchanged for live edges gives a twisted unknot $U(T)$.
\end{lemma}
\pf We need to show that every crossing of $U(T)$ can be undone by a
Reidemeister I move. Given $T$, we can obtain $U(T)$ as follows: first
draw $U$ as if all edges in $G$ are dead; i.e., a regular neighborhood
of $T$, which is a round unknot, up to planar isotopy.
Now for each $e$ not in $T$, we put a crossing there only if
$e$ is live, so this is the only crossing in $cyc(T,e)$, which is a
cycle in a planar graph.
Hence $U$ remains a round unknot after a Reidemeister I move.
Similarly for all live edges $e$ not in $T$.  For all live edges $f$
in $T$, redo this argument for the dual tree $T^*$ using
$cut(T,f)=cyc(T^*,f^*)$.  Therefore, $U(T)$ is isotopic in the plane
to the round unknot after a sequence of Reidemeister I moves.  \done

If $U$ is a partial smoothing of $D$, let $\sigma(U) = \# A$-smoothings $- \#B$-smoothings.
\begin{thm}\label{sum_unknots}
Let $D$ be any connected link diagram, and let $G$ be its Tait graph
with any order on its edges.  There exists a partial skein resolution
tree $\T$, whose leaves are twisted unknots that correspond to
spanning trees of $G$.
If $U$ corresponds to $T$, then $\mu(T)=A^{\sigma(U)}(-A)^{3w(U)}$.
\end{thm}
\pf To construct $\T$, we order crossings of $D$ in the
reverse order to the edges of $G$.  Let the root of $\T$ be $D$.  A
crossing is called {\em nugatory} if either its $A$ or $B$ smoothing
disconnects the diagram.  We smooth the crossings of $D$ in order,
such that at every branch we leave nugatory crossings unsmoothed and
go to the next crossing.  Stop when all subsequent crossings are
nugatory.  Since a diagram is a twisted unknot if and only if all
crossings are nugatory, the leaves of $\T$ are twisted unknots.

From any twisted unknot $U$ in $\T$, we can obtain a
spanning tree $T(U)$ of $G$ by using Table \ref{Table2}, where the
signs below the live edges indicate the writhe of the crossing.
By Lemma \ref{treetounknot}, $U=U(T(U))$.
Each live edge determines the writhe of its crossing in $U$,
hence $\mu(T)=A^{\sigma(U)}(-A)^{3w(U)}$. \done

The activity word for $T$ determines a partial smoothing $U(T)$.
Live edges are not smoothed, denoted below by $*$.

\begin{defn}\label{po}
Let $D$ be any connected link diagram with $n$ ordered crossings.
For any spanning trees $T, T'$ of $G$, let $(x_1,\ldots,x_n)$ and $(y_1,\ldots,y_n)$ be the corresponding partial smoothings of $D$.
We define a relation $T > T'$, or equivalently,
$(x_1,\ldots,x_n) > (y_1,\ldots,y_n)$
if for each $i$, $y_i=A$ implies $x_i = A$ or $*$, and there exists $i$ such that $x_i=A$ and $y_i=B$.
By Proposition \ref{potrans}, the transitive closure of this relation gives a partial order, also denoted by $>$.
We define $P(D)$ to be the poset of spanning trees of $G$ with this partial order.
\end{defn}

\begin{prop}\label{potrans}
If $T >\cdots > T'$ then $T \neq T'$.
\end{prop}
\pf
We can draw $\T$ such that the $A$--smoothing of $c_i$ is $2^{-i}$
units to the left of its parent node, and the $B$--smoothing of $c_i$
is $2^{-i}$ units to the right.  For any $T >\cdots > T'$, $T$ is to
the left of $T'$.
\done

Note that $P(D)$ always has a unique maximal tree and unique minimal
tree, whose partial smoothings contain the all-$A$ and all-$B$
Kauffman states, respectively.

\section{Spanning tree complex}\label{sec3}

Every spanning tree $T$ of a Tait graph $G$ with ordered edges
corresponds to an activity word, which in turn corresponds to a
twisted unknot $U$. 
Let $w(U)$ denote the writhe of $U$, $V(G)$ denote the number of
vertices of $G$ and let $E_{\pm}(G)$ denote the number of positive or
negative edges of $G$.  Given $D$, we require that the checkerboard
coloring be chosen such that $E_+(G)\geq E_-(G)$.

\begin{defn}
  Let $D$ be a connected knot diagram with ordered crossings, and let
  $G$ be its ordered Tait graph.  For any spanning tree $T$ of $G$, we
  define bigradings
$$ u(T)  = -w(U) =   \# L - \# \ell - \# \bar{L} + \# \bar{\ell} \quad {\rm and }\quad
v(T) =  E_+(T)=   \# L + \# D $$
Define $\C(D)=\oplus_{u,v} \C_v^u(D)$, where
$\C_{v}^{u}(D) = \Z\kb{T\subset G |\; u(T)=u,\; v(T)=v }$. \\
Define $\UC(D)=\oplus_{u,v} (\C^u_v(D) + \overline{\C}^{u+2}_{v+1}(D))$,
where $\overline{\C}^{u+2}_{v+1}(D)\cong \C^u_v(D) $.
\end{defn}

\begin{prop}\label{muT}
For any differential $\del: \C_v^u \to \C_{v-1}^{u-1}$, the Jones polynomial
can be expressed as the graded Euler characteristic of $\{\C(D), \del\}$ and of $\{\UC(D), \del\}$:
\begin{eqnarray*}
V_D(t) &=& (-1)^{w(D)}t^{\frac{3w(D)+k}{4}}\ \chi(\C(D)) \\
(t^{1/2}+t^{-1/2})V_D(t) &=& (-1)^{w(D)}t^{\frac{3w(D)+k+2}{4}}\ \chi(\UC(D))
\end{eqnarray*}
where $w(D)$ is the writhe of $D$ and $k= E_+(G) - E_-(G) + 2(V(G)-1)$.
\end{prop}
\pf Let $G$ be the Tait graph of $D$, and let $T$ be any spanning tree of $G$.
By Table \ref{Table2}, the weight of $T$ is given as follows:
\[L^p D^q \ell^r d^s \bar{L}^x \bar{D}^y \bar{\ell}^z \bar{d}^w \quad \Rightarrow \quad \mu(T)= (-1)^{p+r+x+z}A^{-3p+q+3r-s+3x-y-3z+w} \]
Since $T$ is a tree, we have $ p+q+x+y = V(G)-1$ and $r+s+z+w = E(G) - V(G) + 1$.
Also $p+q+r+s = E_+(G)$, $x+y+z+w = E_-(G)$.
Let $k= E_+(G) - E_-(G) + 2(V(G)-1)$.
Since $u = p-r-x+z$, and $v = p+q$,
$\mu(T)= (-1)^u A^{-4(u -v) - k}$.
$$\kb{D} = \sum_{T \subset G} \mu(T) = A^{-k} \sum_u (-1)^u \sum_v A^{-4(u - v)} |\C_v^u|$$
For $t=A^{-4}$, $V_D(t)=(-A)^{-3w(D)}\kb{D}$, so the first result follows.
\begin{eqnarray*}
\chi(\UC(D))&=& \sum_u (-1)^u \sum_v t^{(u-v)}(|\C^u_v|+|\overline{C}^{u+2}_{v+1}|)\\
&=&  \sum_{u,v} (-1)^u t^{(u-v)}|\C^u_v| + t^{-1} \sum_{u,v} (-1)^{u+2} t^{((u+2)-(v+1))}|\overline{C}^{u+2}_{v+1}| \\
&=&(1+t^{-1}) \sum_{u,v} (-1)^u t^{(u-v)}|\C^u_v| \doteq (t^{1/2}+t^{-1/2})V_D(t)
\end{eqnarray*}
The final equality is up to multiplication by
$(-1)^{w(D)}t^{\frac{3w(D)+k+2}{4}}$. \done  

Let $\widetilde{C}(D)=\{\widetilde{C}^{i,j}(D), \del\}$ denote the
reduced Khovanov complex as in \c{Viro}, with
$\widetilde{H}^{i,j}(\bigcirc)=\Z^{(0,-1)}$, where $\bigcirc$ denote
the round unknot.
For chain complexes $X$ and $Y$, $X$ is a {\em deformation retract} of $Y$ if there exist
chain maps $r:Y\to X$ and $f:X\to Y$, such that $r\circ f=id_X$, and a chain homotopy $F:Y\to Y$, such
that $\del_Y F +F\del_Y=id_Y - f\circ r$. Then $r$ is called a retraction.

\begin{thm}\label{mainthm}
For a knot diagram $D$, there exists a spanning tree complex
$\C(D)=\{\C_v^u(D), \del\}$ with $\del$ of bi-degree $(-1,-1)$
that is a deformation retract of ${\widetilde C}(D)$.
In particular, if $w$ is the writhe of $D$, and $k=E_+(G) - E_-(G) + 2(V(G)-1)$,
\begin{equation}\label{uvij}
\widetilde{H}^{i,j}(D;\Z) \cong H_v^u(\C(D);\Z),\quad u = j-i-w+1,\ v = j/2 -i-(w-k-2)/4
\end{equation}
\end{thm}

\begin{thm}\label{unred}
There exists an unreduced spanning tree complex $\UC(D)=(\UC^u_v(D),\del)$
with $\del$ of bi-degree $(-1,-1)$ that is a deformation retract of the (unreduced) Khovanov complex.
In particular, with indices related as in (\ref{uvij}),
  $H^{i,j}(D;\Z) \cong H_v^u(\UC(D);\Z)$.
\end{thm}

\begin{figure}
 \begin{center}
 \psfrag{-}{$-$}
 \psfrag{+}{$+$} \psfrag{=}{$=$}
 \includegraphics[height=0.9in]{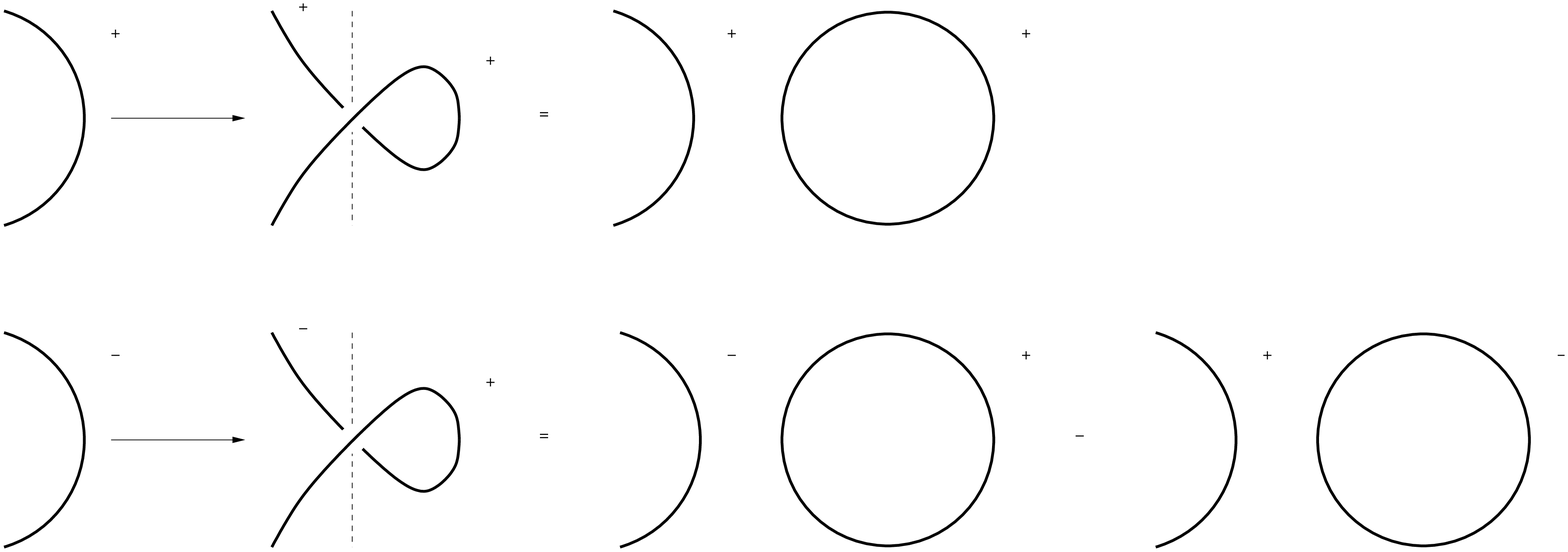}
\hspace*{1.3cm}
 \includegraphics[height=0.9in]{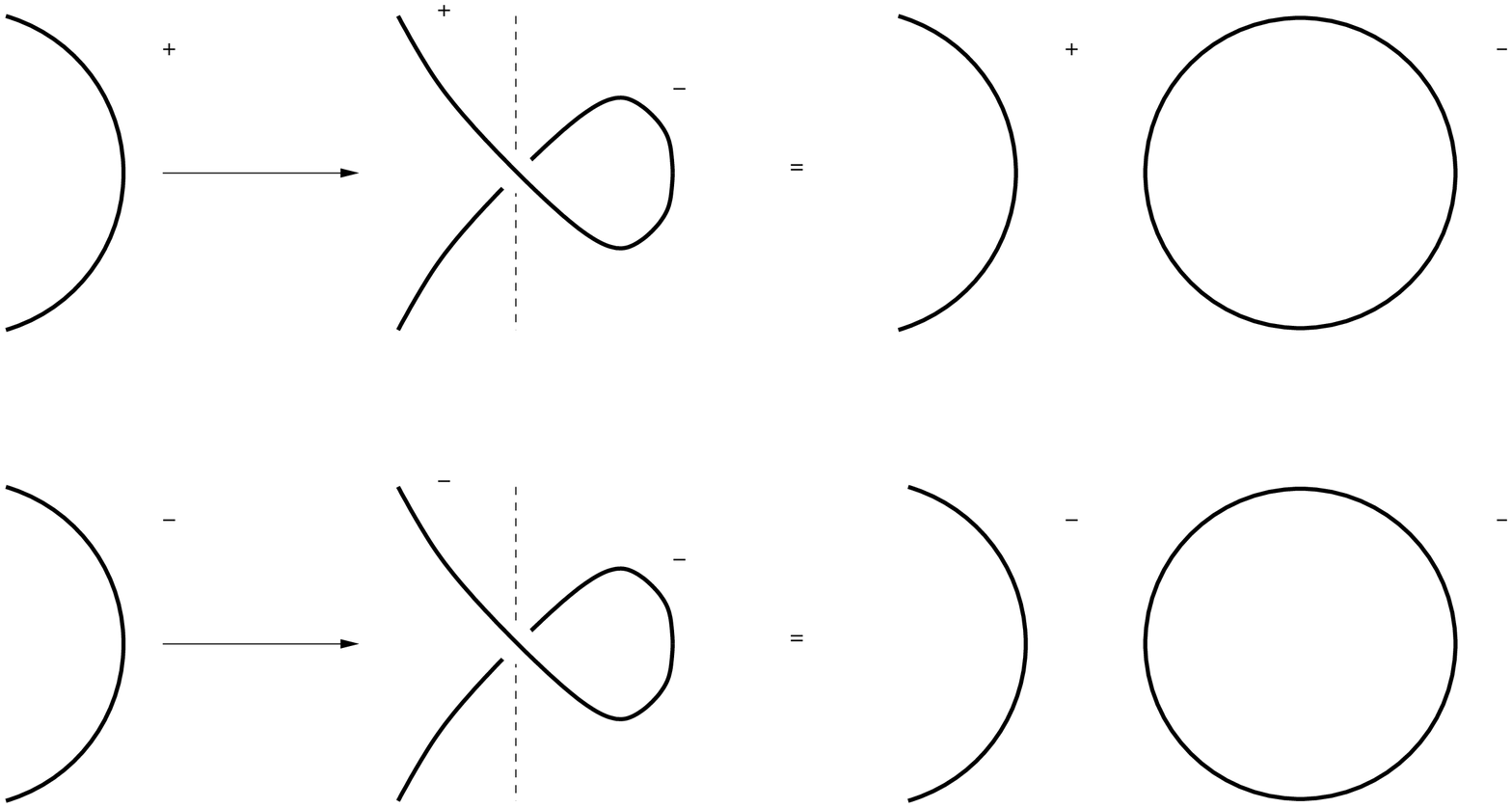}
\end{center}
 \caption{Jacobsson rules for a positive and negative twist}
 \label{postwist}
 \end{figure}

For a twisted unknot $U$, $\widetilde{C}(U)$ is contractible.  Its
homology is generated by a single generator in degree $(i,j)=(0,-1)$,
which is given by iterating the four {\em Jacobsson rules}: Starting
from $\bigcirc$, by a sequence of positive and negative twists, we
obtain $U$, and Figure \ref{postwist} indicates how to change the
enhanced state for each twist, starting with the round unknot enhanced
by a $+$ sign, $\bigcirc^+$, which generates 
${\widetilde C}(\bigcirc)\cong \Z^{(0,-1)}$.

\begin{defn}\label{Jrules}
For any twisted unknot $U$, we define its fundamental cycle $Z_U \in
\widetilde{C}(U)$ to be the linear combination of maximally
disconnected enhanced states of $U$ given by the Jacobsson rules. Let
$f_U: {\widetilde C}(\bigcirc) \to {\widetilde C}(U)$ be defined by
$f_U(\bigcirc^{+})=Z_U$.
\end{defn}

Let $w(U)$ be the writhe of $U$, let
$\sigma$ be the difference of $A$-smoothings and $B$-smoothings, and
let $\tau$ be the difference of positive and negative enhancements, as
in \c{Viro}.  By Figure \ref{postwist}, 
the Jacobsson rules have the following effect for each added twist:
\begin{center}
\begin{tabular}{lccc}
Positive twist: & $w\to w+1$ & $\sigma\to \sigma +1$ & $\tau\to\tau +1$\\
Negative twist: & $w\to w-1$ & $\sigma\to \sigma -1$ & $\tau\to\tau -1$
\end{tabular}
\end{center}
The grading for any enhanced state in
$\widetilde{C}^{i,j}(U)$ is given by $i=(w-\sigma)/2$ and
$j=i+w-\tau$, which are preserved under the Jacobsson
rules.  By Lemma \ref{collapse2} below, $f_U$ is a grading-preserving chain homotopy.
Let $\iota : {\widetilde C}(U) \to {\widetilde C}(D)$ be the inclusion
of enhanced states of $U$ into enhanced states of $D$ given by the
grading shifts $\iota(s^{i,j})=s^{i',j'},$ where
$i'=i+\frac{w(D)-w(U)-\sigma(U)}{2},\
j'=j+\frac{3(w(D)-w(U))-\sigma(U)}{2}$.
For any spanning tree $T_k$, we define $\U_k = \iota({\widetilde
C}(U(T_k))) \subset {\widetilde C}(D)$.

{\em Proof of Theorem \ref{mainthm}:}\quad
Fix any order on the crossings of $D$ and a basepoint on $D$ away from the crossings.
We embed $\C(D)$ into ${\widetilde C}(D)$ as bigraded groups.  For
each generator $T\in \C_v^u(D)$, let $U=U(T)$.
Let $\phi: \C(D)\to {\widetilde C}(D)$ be defined by 
$\phi(T) = \iota\left(Z_{U(T)}\right)$.

For given $(u(T),v(T))$, we compute the $(i,j)$ degree of $\phi(T)$ in $\widetilde{C}(D)$.
Let $\sigma$ and $\sigma_U$ denote the signature $(\# A - \# B)$ for an enhanced state in $\widetilde{C}(D)$ and $\widetilde{C}(U)$, respectively.
Since $Z_{U(T)}\in\widetilde{C}^{0,-1}(U)$,
we have $w(U)-\sigma_U=0$ and $w(U)-\tau=-1$.
Since $u(T)=-w(U)$, we have $\sigma_U = -u(T)$ and $\tau = 1 - u(T)$.
Since $U$ is obtained from $D$ by smoothing all crossings on dead edges,
using the notation of Proposition \ref{muT}, $\sigma = \sigma_U + (q-s-y+w)$.
Since $q-s-y+w = -u +4v -k(D)$, $\tau$ and $\sigma$ of $Z_{U(T)}$ are 
$\tau = 1 - u(T)$ and $\sigma = -2u(T) + 4v(T) - k(D)$.
Therefore, $\phi(T)$ has the following $(i,j)$ degree in $\widetilde{C}(D)$:  If $w=w(D)$ and $k=k(D)$,
\begin{equation}\label{ijuv}
 i= \frac{w-\sigma}{2} = u-2v + \frac{w+k}{2} \quad {\rm and}\quad j= i+w-\tau= 2u-2v + \frac{3w+k-2}{2}
\end{equation}
Solving for $u$ and $v$, we obtain (\ref{uvij}).

We now order the spanning trees of $G$ as $T_k,\ 1 \leq k \leq s$, such that if $T_{k_1} > T_{k_2}$ then $k_1 > k_2$.
We proceed by a sequence of elementary collapses of each unknot's complex to its fundamental cycle starting from the minimal tree.
Lemma \ref{collapse3} provides the crucial fact that any elementary collapse in $\U_k$ does not change
incidence numbers in $\U_c$ for any $c\neq k$.
This fact permits us to repeatedly apply Lemma \ref{collapse2}:
Starting with $C_0 = {\widetilde C}(D)$,
we get a sequence of complexes $C_k,\, 0 \leq k \leq s$, and retractions $r_k: C_0 \rightarrow C_{k}$.
Each $C_{k+1}$ is obtained from $C_k$ by a sequence of elementary collapses by deleting all generators in
$r_k(\U_k)$ except the fundamental cycle.
Finally, $C_s$ is a complex whose generators are in $1-1$ correspondence with the spanning trees of $G$ and
${\widetilde H}^{i,j}(C_s) \cong {\widetilde H}^{i,j}(D)$.

The map $r_s \circ \phi: \C_v^u(D) \rightarrow C_s$ is a graded
isomorphism of groups, with the indices related by (\ref{uvij}).  The
induced differential on the spanning tree complex $\C_v^u(D)$ now gives
$H_v^u(\C(D)) \cong \widetilde{H}^{i,j}(D)$ with the indices related by
(\ref{uvij}).  In the version of Khovanov homology in \c{Viro}, the
differential on ${\widetilde C}(D)$ has bi-degree $(1,0)$, so by
(\ref{uvij}) the differential on $\C_v^u(D)$ has bi-degree $(-1,-1)$.
The retraction from the reduced Khovanov complex to the spanning tree complex is given by
\begin{equation}\label{retract}
r = \left(r_s\circ\phi\right)^{-1} \circ r_s: \widetilde{C}^{i,j}(D) \to \C_v^u(D)
\end{equation}
where $r(\U_k)=T_k,\ r\circ\phi =id$, and the indices are related by (\ref{uvij}).
\done

For a complex  $(C,\del)$ over $\Z$,
we say $x$ is \emph{incident} to $y$ in $(C,\del)$ if
$\kb{\del x, y}\neq 0$ and their \emph{incidence number} is $\kb{\del x, y}$.

\begin{lemma}\label{podel}
The differential $\del$ on ${\widetilde C}(D)$ respects the partial
order in Definition \ref{po}: 
$(i.)$ If $\kb{\del x, y}\neq 0$ for any $x\in\U_1$ and $y\in\U_2$, then $T_1 > T_2$.
$(ii.)$ If $T_1$ and $T_2$ are not comparable or $T_2 > T_1$, then $\kb{\del x, y}=0$ for all $x\in\U_1$ and $y\in\U_2$.
\end{lemma}
\pf If $\kb{\del x, y}\neq 0$ then any partial smoothing that contains
these states satisfies $(x_1,\ldots,x_n) > (y_1,\ldots,y_n)$ as in
Definition \ref{po}.  \done

\begin{lemma}\label{collapse1} {\bf (Elementary collapse)}
Let $(C,\del)$ be a finitely generated chain complex over $\Z$.  Let $x,\ y$
be generators, such that $x\in C_k,\ y\in C_{k-1}$ with $\kb{\del x, y}=\pm 1$.
Then there exists a complex $(C',\del)$, such that
$C'$ is generated by the same generators as $C$ except for $x,\, y$,
and there is a retraction $r:C \rightarrow C'$.
\end{lemma}
\pf
Fix bases $E_n$ of $C_n$ with
$E_{k-1}=\{y,\, y_1,\ldots,y_{d_{k-1}}\}$, $E_k=\{x,\, x_1,\ldots,x_{d_k}\}$.
For $n\geq 0$, define $r_n:C_n \rightarrow C_n$ as follows:
For any $v\in C_n$, $r_n(v)=v$ if $n\neq k,k-1$,
$$ r_n v = v - \frac{\kb{v,y}}{\kb{\del x,y}}\del x \text{ if } n = k-1
\quad {\rm and }\quad
r_n v = v - \frac{\kb{\del v,y}}{\kb{\del x,y}}x \text{ if } n = k $$
Define $r:C \rightarrow C$ as $r|_{C_n}=r_n$. Then $r$ is a chain map and hence
its image is a subcomplex. Let $(C',\del)=(r(C),\del)$ .
If $\del x = \lambda y + Y$, with $\lambda=\pm1$ and $\kb{y,Y}=0$,
then for $i\geq 1,\ r_{k-1}(y_i)=y_i$, and $r_{k-1}(y)= -\lambda Y$.
For $i\geq 1,\ r_k(x_i) = x_i - \lambda\kb{\del x_i,y} x$, and $r_k(x) = 0$.
It follows that $r:C \rightarrow C'$ is a retraction. 
\done

\begin{lemma}\label{collapse2}
Let $U$ be a twisted unknot. There exists a sequence of elementary collapses $r_U: {\widetilde C}(U)\to {\widetilde C}(\bigcirc)$, 
such that $r_U\circ f_U=id$ and $f_U\circ r_U\simeq id$. 
\end{lemma}
\pf In essence, this result follows from invariance of Khovanov
homology under the first Reidemeister move \c{Khovanov}, but we
explicitly provide the elementary collapses.
The proof is by induction on the number of crossings of $U$.  Suppose
$U'$ is obtained from $U$ by adding one kink, hence one crossing $c$.
The markers below refer to $c$, and the signs to the enhancements near
$c$ in the order they appear in Figure \ref{postwist}. 

{\bf Positive twist} 
The $A$-smoothing of $U'$ at $c$ results in a new loop; the $B$-smoothing does not.
For every enhanced state $v^+$ of $U$, collapse the pair $A^{+-} {\to} B^+$.
For every enhanced state $v^-$ of $U$, collapse the pair $A^{--} {\to} B^-$.
By Lemma \ref{collapse1}, $r(A^{++})=v^+$ and since $A^{-+} {\to} B^+$, we get $r(A^{-+}-A^{+-})=v^-$.

{\bf Negative twist} 
The $B$-smoothing of $U'$ at $c$ results in a new loop; the $A$-smoothing does not.
For every enhanced state $v^+$ of $U$, collapse the pair $A^+ {\to} B^{++}$.
For every enhanced state $v^-$ of $U$, collapse the pair $A^- {\to} B^{-+}$.
By Lemma \ref{collapse1}, $r(B^{+-})=v^+$ and $r(B^{--})=v^-$.

Let $f:{\widetilde C}(U)\to {\widetilde C}(U')$ be the following map:
For any $s\in{\widetilde C}(U)$, let $f(s)$ be the linear combination of states given by
Figure \ref{postwist}.  From the change in $w,\ \sigma$ and $\tau$, $f$ is grading-preserving.
Moreover, $f$ is an iterated Jacobsson map: $f \circ f_U=f_{U'}$.
The elementary collapses above show that
$ {\widetilde C}(U) \overset{f}{\rightarrow}  {\widetilde C}(U') \overset{r}{\rightarrow} {\widetilde C}(U) $
with $r\circ f=id$ and $f\circ r \simeq id$.
Starting with $U=\bigcirc$, the result follows by induction. \done

\begin{lemma}\label{collapse3}
Let $D$ be a connected link diagram and let $G$ be its Tait graph. Let
$T_1$ and $T_2$ be distinct spanning trees of $G$.  Then in
${\widetilde C}(D)$, any elementary collapse as in Lemma
\ref{collapse1} of $x_1,y_1 \in \U_1$ will not change the incidence
number between any $x_2,y_2 \in \U_2$.
\end{lemma}
\pf As in the proof of Lemma \ref{collapse1}, $\kb{\del x_1,y_1} = \lambda \in\{\pm1\}$.  The image
of $x_2$ after the elementary collapse of $x_1, y_1$ is $x_2' = r(x_2)
= x_2 - \lambda\kb{\del x_2,y_1}x_1$.  Hence, 
$\kb{\del x_2',y_2}=\kb{\del x_2,y_2} - \lambda\kb{\del
x_2,y_1}\kb{\del x_1,y_2}$.
By Lemma \ref{podel}, if $T_1>T_2$ then $\kb{\del x_2,y_1}=0$, and otherwise $\kb{\del x_1,y_2}=0$.
Thus, $\kb{\del x_2',y_2}=\kb{\del x_2,y_2}$. \done

{\em Proof of Theorem \ref{unred}:}\quad
For (unreduced) Khovanov homology, $H^{i,j}(\bigcirc;\Z)=\Z^{0,1}\oplus\Z^{0,-1}$.
So the Khovanov complex for any twisted unknot $U$ is chain homotopic
to the complex with two generators in degrees $(i,j)=(0,\pm
1)$.  Hence, for every $T$, there are two fundamental
cycles for $U(T)$, and two corresponding generators:
$T_+$ in grading $(u(T),v(T))$ and $T_-$ in grading $(u(T)+2,v(T)+1)$.
Lemmas \ref{podel}, \ref{collapse2} and \ref{collapse3} now extend to the unreduced Khovanov complex, and
the rest of the proof follows from that of Theorem \ref{mainthm}.
\done

\section{Spanning tree filtration and spectral sequence}\label{sec5}

The poset of spanning trees $P$ given in Definition \ref{po}, together with
Proposition \ref{potrans} and Lemma \ref{podel}, provide a partially ordered filtration of ${\widetilde C}(D)$ indexed by $P$:
Let $\psi: P\to {\widetilde C}(D)$ be defined by
$\psi(T) = +_{T \geq T_i} \U_i$.
This defines a decreasing linearly ordered filtration on ${\widetilde C}(D)$ as follows.
Let $\{S_j\}$ be the set of maximal descending ordered sequences of spanning trees in $P$, and let $T^j_k$ denote the $k$-th element of $S_j$, so that for all $j$, $T^j_1$ is the unique maximal tree in $P$.
Define
$F^p {\widetilde C}(D) = +_j\ \psi(T^j_p)$.

\begin{thm}\label{SS}
For any knot diagram $D$, there is a spectral sequence $E_r^{*,*}$
that converges to the reduced Khovanov homology
$\widetilde{H}^{*,*}(D;\Z)$, such that as groups $E^{*,*}_1 \cong
\C^*_*(D)$, and the spectral sequence collapses for $r\leq c(D)$, where
$c(D)$ is the number of crossings.
\end{thm}
\pf
By Lemma \ref{podel}, the differential on ${\widetilde C}(D)$ respects the filtration $\{F^p {\widetilde C}(D)\}$;
i.e., $\partial F^p \subseteq F^p$.
Hence this filtration determines a spectral sequence $\{E^{p,q}_r,\,d_r\}$
that converges to the reduced Khovanov homology.
The associated graded module consists of submodules of ${\widetilde C}(D)$ in bijection with spanning trees:
\begin{equation} \label{E0}
 E^{p,*}_0 = F^p{\widetilde C}(D)/F^{p+1}{\widetilde C}(D) = \oplus_k\ \U_k
\end{equation}
It follows from the filtration that for any $p$, if
$\U_1, \U_2\subset E^{p,*}_0$, then $T_1$ and $T_2$ are not comparable in $P$.  Hence,
by Lemma \ref{podel}($ii$), $d_0: E^{p,q}_0 \to E^{p,q+1}_0$ satisfies $d_0(\U_k)\subset\U_k$ for every $k$.
This implies that (\ref{E0}) is a direct sum of complexes $\U_k$.
By Lemma \ref{collapse2}, each complex $\U_k$ has homology generated by $\phi(T_k)$.
Therefore, $E_1$ is isomorphic as a group to the spanning tree complex:
$$E^{*,*}_1 =H^*(F^p/F^{p+1},d_0)= \oplus_k\ H^*(\U_k) \cong \C^*_*(D)$$
Since the length of any ordered sequence in $P$ is at most the number of crossings $c(D)$,
it follows that the spectral sequence collapses for $r\leq c(D)$.
\done

For field coefficients, Theorem \ref{SS} provides another proof that a differential exists on
the spanning tree complex $\C(D)$ that makes it chain homotopic to $\widetilde{C}(D)$:

\begin{thm}\label{mainthmoverfield} For coefficients in a field $\mathbb{F}$, there exists a differential on $\C(D)$
such that $\widetilde{H}^{*,*}(D;\mathbb{F}) \cong H^*_*(\C(D);\mathbb{F})$.
\end{thm}
\pf By Lemma 4.5 in \cite{Rasmussenthesis}, there exists a unique filtered complex $C'$ which is chain homotopic to
$\widetilde{C}(D)$ and $C' \cong H^*(F^p\widetilde{C}(D)/F^{p+1}\widetilde{C}(D))$.
Theorem \ref{SS} implies that $C' \cong H^*(E_0^{*,*}) \cong E_1^{*,*} \cong \C(D)$.
\done

\section{Alternating and almost alternating links}\label{sec7}
We give a simple new proof for theorems of Lee \c{Lee} and Shumakovitch \c{Shurik} for the reduced Khovanov homology
using the spanning tree complex.
\begin{thm} \label{altredKH} The reduced Khovanov homology of an alternating knot is determined by its
 Jones polynomial and signature, and it has no torsion.
\end{thm}
\pf An alternating diagram $D$ can be checkerboard colored so that its Tait graph $G$ has all positive edges.
For any spanning tree $T$ of $G$, $v(T)=E_+(T)=E(T)=V(G)-1$.
Since the $v$-grading is constant for all spanning trees,
all the generators in the spanning tree complex $\C(D)$ are in one row. Since the differential on $\C(D)$
has degree $(-1,-1)$, it is trivial. Hence by Theorem \ref{mainthm},
$\widetilde{H}^{i,j}(D;\Z) \cong H_v^u(\C(D);\Z) \cong \C_v^u(D)$.
Therefore, the homology has no torsion.  The Betti numbers are determined by the Jones polynomial:
If $c(D)$ is the number of crossings of $D$, Proposition \ref{muT} implies that $|\C_v^u(D)|=a_{u-v+\frac{3w(D)+c(D)+2v}{4}}$, where $V_D(t)=\sum a_n t^n$, and we use that $k(D)=E(G)+2(V(G)-1)=c(D)+2v$.
By \c{Traczyk}, the signature of the knot $\sigma(D)=\frac{c(D)-w(D)}{2}-|s_B|+1$, where $s_B(D)$ is the Kauffman state with all $B$ markers.
Since $D$ is alternating, $|s_B(D)|=V(G)=v+1$.  Therefore,
$v = \frac{c(D)-w(D)}{2} - \sigma(D)$.
\done

\begin{rmk}\label{ijrmk}
Using (\ref{ijuv}), the above proof implies that for non-split alternating links,
$$ j-2i= 2(V(G)-1) +\frac{w(D)-k(D)}{2}-1 = v - \frac{c(D)-w(D)}{2} - 1 = -\sigma(D)-1 $$
\end{rmk}

A link is {\em $k$-almost alternating} if it has a non-nugatory
diagram which is alternating after $k$ crossing changes, and does
not have one after $k-1$ crossing changes.
The bigraded homology of a link is {\em $k$-thick} if the
nontrivial homology groups lie on $k$ adjacent lines.  We give 
a simple new proof for theorems
about the support of Khovanov homology for alternating and $k$-almost
alternating links obtained respectively by Lee \c{Lee} and Asaeda,
Przytycki \c{AP}.
We proved a more general result in terms of ribbon graph genus in \c{dkh}.
Another proof also appeared in Manturov \c{Manturov}.

\begin{thm}\label{thick}
$(i)$ The Khovanov homology of a non-split alternating link $L$ is at most $2$-thick, and lies on the lines
$j-2i=-\sigma(L)\pm 1$. Its torsion lies on the line $j-2i=-\sigma(D)- 1$.
$(ii)$ The Khovanov homology of a non-split $k$-almost alternating link $L$
  is at most $(k+2)$-thick, and its reduced Khovanov homology
  is at most $(k+1)$-thick.
\end{thm}

\pf $(i)\ $ For an alternating diagram $D$, $\UC(D)$ lies on
two lines, $v=V(G)-1$ and $v=V(G)$.  From Remark \ref{ijrmk}, the
homology lies on the lines $j-2i=-\sigma(D)\pm 1$.  Moreover, since
the differential on $\UC(D)$ has degree $(-1,-1)$, any
torsion in the homology must lie on the line $j-2i=-\sigma(D)- 1$.

$(ii)\ $ A $k$-almost alternating link or its mirror image has a Tait graph $G$ with exactly $k$ negative edges, such
that $k\leq E(G)/2$. For any spanning tree $T$ of $G$, $v(T) = E_+(T)$, so $v(T)$ takes at most $(k+1)$ values.
Since $\UC(D)$ has at most $(k+2)$ rows, 
$H_v^u(\UC(D))$ has at most $(k+2)$ rows.
The result now follows from Theorems \ref{mainthm} and \ref{unred}.
\done

\section{Example}\label{example}

\begin{tabular}{{m{2in} m{3in}}}
As an example we use a 4-crossing diagram of the trefoil. Here is the diagram $D$ and its Tait graph $G$.
Below we show all the spanning trees of $G$  with their activity words, $(u,v)$-gradings and partial smoothings.
&
\makebox[3in]{
\psfrag{1}{\tiny{$1$}}
\psfrag{2}{\tiny{$2$}}
\psfrag{3}{\tiny{$3$}}
\psfrag{4}{\tiny{$4$}}
\psfrag{3b}{\tiny{$\overline{3}$}}
\psfrag{4b}{\tiny{$\overline{4}$}}
 \includegraphics[height=1in]{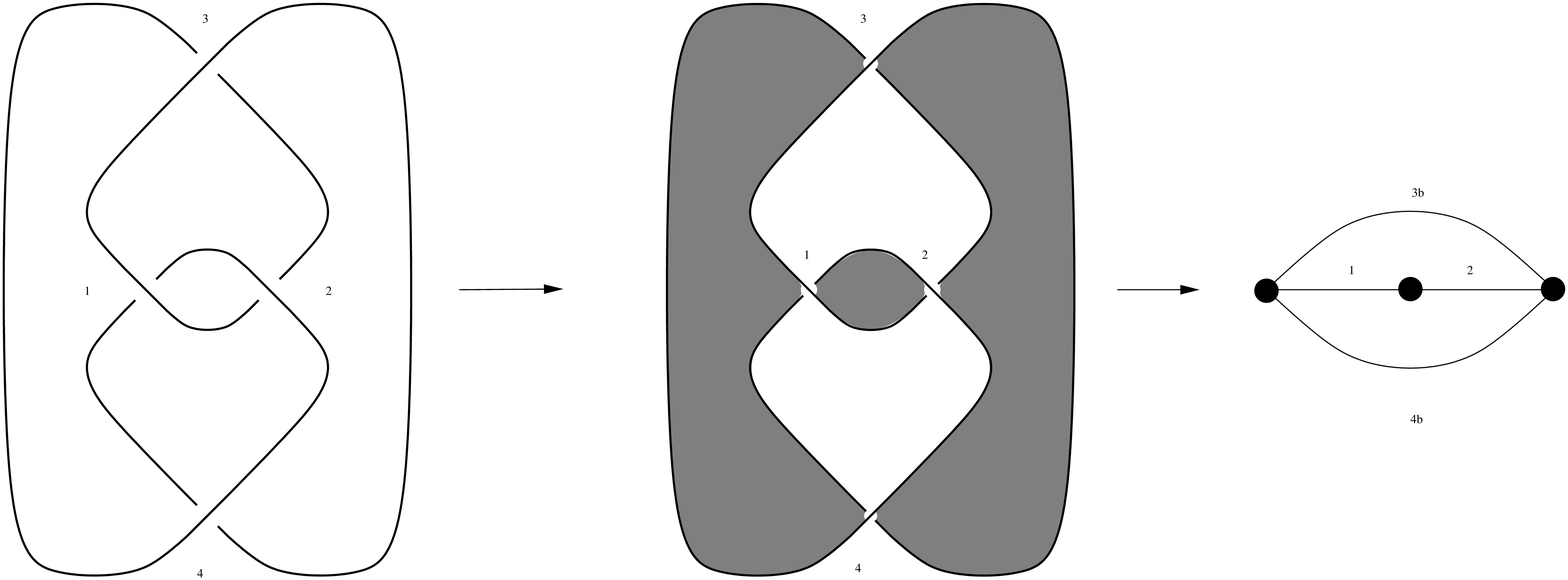}
}
\end{tabular}

\begin{center}
\psfrag{1}{\tiny{$1$}}
\psfrag{2}{\tiny{$2$}}
\psfrag{3}{\tiny{$3$}}
\psfrag{4}{\tiny{$4$}}
\psfrag{3b}{\tiny{$\overline{3}$}}
\psfrag{4b}{\tiny{$\overline{4}$}}
\psfrag{t1}{\footnotesize{$T_1=\ell D \overline{D d}$}}
\psfrag{t2}{\footnotesize{$T_2=\ell D \overline{\ell D}$}}
\psfrag{t3}{\footnotesize{$T_3=L d \overline{D d}$}}
\psfrag{t4}{\footnotesize{$T_4=L d \overline{\ell D}$}}
\psfrag{t5}{\footnotesize{$T_5=L L \overline{d d}$}}
\psfrag{U1}{\footnotesize{$(-1,1)$}}
\psfrag{U2}{\footnotesize{$(0,1)$}}
\psfrag{U3}{\footnotesize{$(1,1)$}}
\psfrag{U4}{\footnotesize{$(2,1)$}}
\psfrag{U5}{\footnotesize{$(2,2)$}}
\psfrag{P1}{\footnotesize{$*ABA$}}
\psfrag{P2}{\footnotesize{$*A*B$}}
\psfrag{P3}{\footnotesize{$*BBA$}}
\psfrag{P4}{\footnotesize{$*B*B$}}
\psfrag{P5}{\footnotesize{$**AA$}}
\includegraphics[width=4.25in]{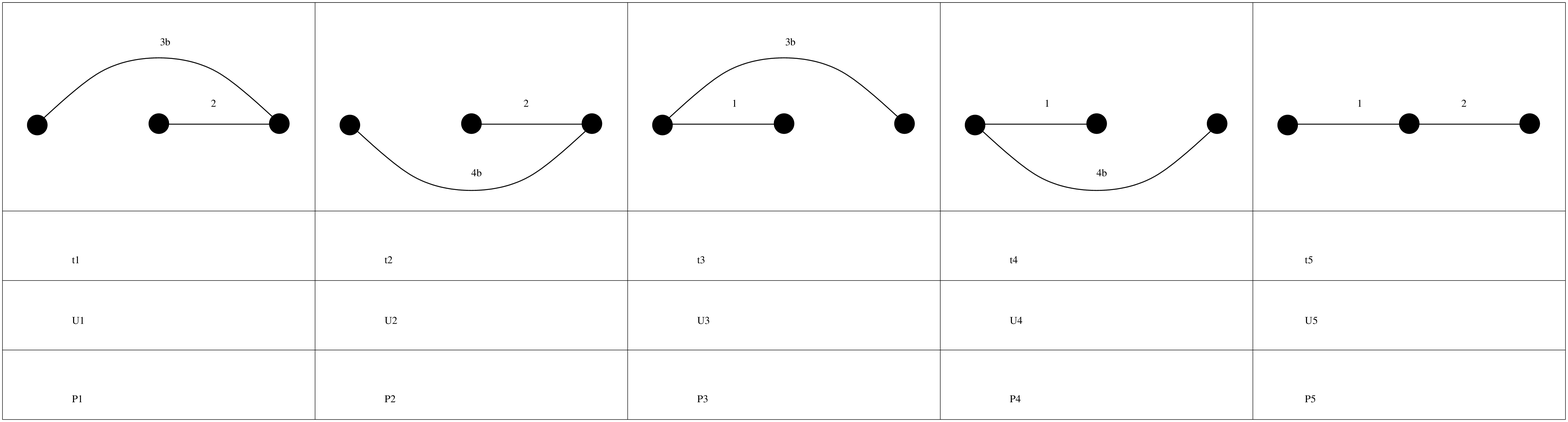}
 \end{center}
There are two maximal sequences in the partial order: $T_5>T_1>T_2>T_4$ and $T_5>T_1>T_3>T_4$.
The associated graded module $E^{p,q}_0$ is: $E^{1,q}_0=\widetilde{U}_5$, $E^{2,q}_0=\widetilde{U}_1$,
$E^{3,q}_0=\widetilde{U}_2\oplus \widetilde{U}_3$, $E^{4,q}_0=\widetilde{U}_4$, with $q$ determined
by $p+q=i=u-2v+2$. We show the $E_1$, $E_2$ and $E_3$ pages of the spectral sequence, which collapses at $E_3$.
The generators and differentials are shown with $(u,v)$-gradings, which determine 
the $(p,q)$-gradings for the spectral sequence. 

\begin{center}
\psfrag{t1}{\footnotesize{$T_1$}}
\psfrag{t2}{\footnotesize{$T_2$}}
\psfrag{t3}{\footnotesize{$T_3$}}
\psfrag{t4}{\footnotesize{$T_4$}}
\psfrag{t5}{\footnotesize{$T_5$}}
\psfrag{e1}{\footnotesize{$E_1^{*,*}$}}
\psfrag{e2}{\footnotesize{$E_2^{*,*}$}}
\psfrag{e3}{\footnotesize{$E_3^{*,*}$}}
\includegraphics[width=4in]{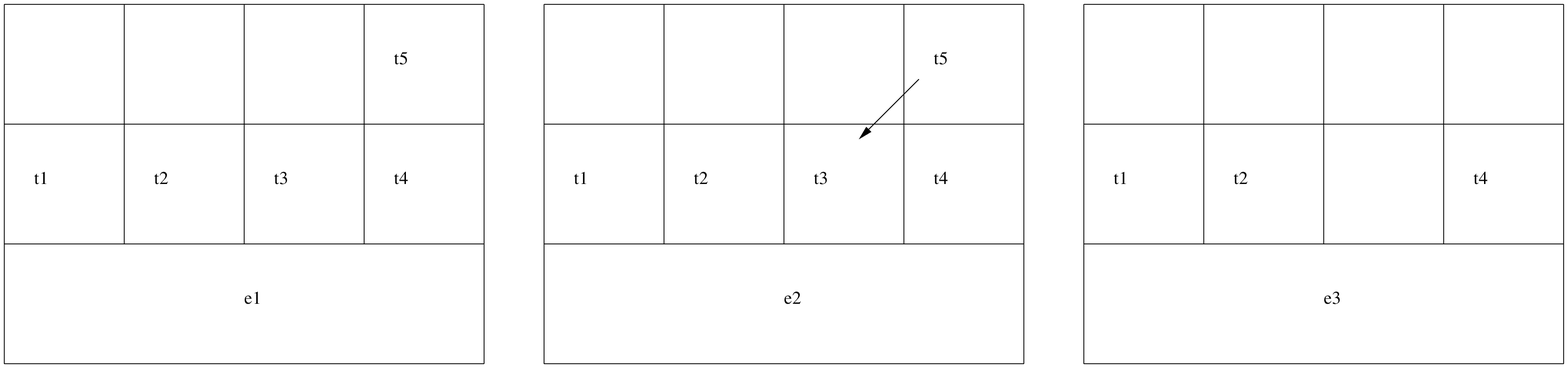}
 \end{center}

\bibliography{resplan}
\bibliographystyle{plain}
\end{document}